\input amstex
\documentstyle{amsppt}
\NoBlackBoxes

\mathchardef\t="0074 \mathchardef\o="006F

\mathchardef\p="0070 \mathchardef\r="0072 \mathchardef\t="0074
\mathchardef\b="0062 \mathchardef\d="0064 \mathchardef\i="0069
\mathchardef\n="006E \mathchardef\w="0077 \mathchardef\s="0073
\def\bt{\mathop\boxtimes\limits}
\def\dis{\displaystyle}

\topmatter
\title On the classification of inductive limits \\ of II$_1$
factors  with spectral gap
\endtitle
\author SORIN POPA \endauthor

\rightheadtext{Classification of inductive limits}

\affil University of California, Los Angeles\endaffil

\address Math.Dept., UCLA, LA, CA 90095-155505\endaddress
\email popa\@math.ucla.edu\endemail

\thanks Supported in part by NSF Grant 0601082.\endthanks

\abstract We consider II$_1$ factors $M$ which can be realized as
inductive limits of subfactors, $N_n \nearrow M$, having spectral
gap in $M$ and satisfying the bi-commutant condition $(N_n'\cap
M)'\cap M=N_n$. Examples are the enveloping algebras associated to
non-Gamma subfactors of finite depth, as well as certain crossed
products of McDuff factors by amenable groups. We use
deformation/rigidity techniques to obtain classification results for
such factors.
\endabstract

\endtopmatter

\document

\heading 1. Introduction \endheading

A von Neumann subalgebra $Q$ of a II$_1$ factor $M$ has {\it
spectral gap in} $M$ if any element $x\in M$ that almost commutes
with $Q$ must be close (in the Hilbert norm given by the trace) to
an element that actually commutes with $Q$. This condition holds
true, for instance, whenever $Q$ has the property (T). But it is
also satisfied under much weaker conditions on $Q$, e.g. if $Q$ does
not have amenable direct summand, provided this is ``compensated''
by certain requirements on how $Q$ sits inside $M$.

In [P06a, P06b], spectral gap was used in combination with
malleability properties of the ambient factor $M$, to derive
rigidity results, using the deformation-rigidity techniques
developed in ([P01, P03, P04a]; see [P06c] for a survey).
In this paper we obtain several more
applications of spectral gap rigidity, by combining it instead with
inductive limits deformation.

Our starting point is the observation that if $Q$ has spectral gap
in $M$ and $M$ is an inductive limit of subfactors, $N_n \nearrow
M$, then $N_n'\cap M$ is approximately contained in $Q'\cap M$, for
large $n$. Furthermore, if $N_n$ and $Q$ satisfy the bicommutant
condition in $M$, i.e. $(Q'\cap M)'\cap M=Q$ and $(N_n'\cap M)'\cap
M=N_n$, then $Q$ follows almost contained in $N_n$ as well. In
particular, by using the {\it intertwining subalgebras} techniques
and notations in [P03], one gets $Q\prec_M N_n$, i.e. $Q$ can be
conjugated into $N_n$ by a unitary element of $M$ (roughly). This
shows that in order for the deformation-rigidity arguments to work,
one needs the II$_1$ factors $M$ to be inductive limits of
subfactors, $N_n \nearrow M$, having spectral gap in $M$ and
satisfying the bi-commutant condition $(N_n'\cap M)'\cap M=N_n$. As
it turns out, the additional condition $[N_n:N_m] < \infty$,
$\forall n>m$, is also needed.

An important class of factors which arise as inductive limits of
subfactors satisfying all these conditions, are the enveloping
algebras $N_\infty$ associated to non-Gamma subfactors of finite
depth $N_{-1}\subset N_0$, defined by $N_\infty = \overline{\cup_n
N_n}$, where $N_n \nearrow N_\infty$ denotes the Jones tower for
$N_{-1}\subset N_0$. The above argument then shows that any
isomorphism between such enveloping factors, is implemented by a
{\it weak equivalence} of the finite depth subfactors they come from
(see Definition 3.1.4$^\circ$ and Corollary 3.6). In other words,
the enveloping algebra $N_\infty$ of a non-Gamma subfactor with
finite depth $N_{-1}\subset N_0$, ``roughly'' remembers
$N_{-1}\subset N_0$.

Another class of factors satisfying the above conditions are the
crossed product factors of the form $T \rtimes \Gamma$, where $T$ is
a {\it s-McDuff factor} (i.e., a factor of the form
$N\overline{\otimes} R$, with $N$ a non-Gamma II$_1$ factor and $R$
the hyperfinite II$_1$ factor), $\Gamma$ is a countable amenable
group and $\Gamma \curvearrowright T$ is a {\it proper} action of
$\Gamma$ on $T$, where by definition this means $\Gamma$ acts
outerly on both $N$ and $R$ (the unique, by [P06], non-Gamma and
hyperfinite components of $T$) and its image in Out$(T)$ is closed.

A concrete class of such examples is obtained as follows. Let $N$ be
a non-Gamma II$_1$ factor (e.g. $N=L(\Bbb F_n)$, for some $2 \leq n
\leq \infty$), $\Gamma$ a countable group and $\theta : \Gamma
\hookrightarrow \text{\rm Out}(N)$, $\rho: \Gamma \hookrightarrow
\text{\rm Out}(R)$ faithful $\Gamma$-kernels, such that their
corresponding H$^3(\Gamma, \Bbb T)$-obstructions satisfy
Ob$(\theta)=\overline{\text{\rm Ob}(\rho)}$ and such that the image
of $\theta(\Gamma)$ in Out$(N)$ is closed. Then
$\sigma=\theta\otimes \rho$ implements a free cocycle action of
$\Gamma$ on $T=N\overline{\otimes} R$ and all the above conditions
are satisfied. Note that by [P91], [PS00], if $N=L(\Bbb F_\infty)$
then, given any countable amenable group $\Gamma$ and any element
$\alpha \in \text{\rm H}^3(\Gamma, \Bbb T)$, there exists $\theta:
\Gamma \hookrightarrow \text{\rm Out}(L(\Bbb F_\infty))$ such that
Ob$(\theta)=\alpha$ and such that $\theta(\Gamma)$ is closed in
$\text{\rm Out}(L(\Bbb F_\infty))$. Also, by [J80], there exists a
faithful $\Gamma$-kernel $\rho:\Gamma \hookrightarrow \text{\rm
Out}(R)$ such that Ob$(\rho)=\overline{\alpha}$, which by [Oc83] is
in fact unique, up to cocycle conjugacy.

Our results show that any isomorphism $\Delta:M_0\simeq M_1$ between
crossed product factors $M_i=T_i \rtimes_{\sigma_i} \Gamma_i$ in
this class, ``virtually'' takes the s-McDuff factors $T_0, T_1$ onto
each other, and cocycle conjugates $\sigma_0, \sigma_1$, as well as
their restrictions $\theta_0, \theta_1$ to the non-Gamma parts of
$T_0, T_1$. Moreover, if $\Gamma_i$ are torsion free, then $\Delta$
comes from an actual cocycle conjugacy of $\theta_0, \theta_1$. In
particular, the group $\Gamma_i$ and the element $\alpha_i=\text{\rm
Ob}(\theta)\in \text{\rm H}^3(\Gamma_i,\Bbb T)$, are isomorphism
invariants for $M_i$. For instance, since H$^3(\Bbb Z^3,\Bbb T)=\Bbb
Z$, if we take $\Gamma=\Bbb Z^3$ and $\alpha \in \Bbb T\setminus
\{\pm 1\}$, then the corresponding factor $M=L(\Bbb F_\infty)
\overline{\otimes} R \rtimes_\sigma \Bbb Z^3$ satisfies $M\not\simeq
M^{\o\p}$ and $M^{\otimes n}$ are non-isomorphic, for $n=1,2,3...$.
This construction should be compared with Connes original examples
satisfying such properties, obtained using his $\chi(M)$ invariant
(see [C75]). We mention that in a paper that we circulate in
parallel ([P09]), we use similar deformation by inductive limits and
spectral gap rigidity to calculate $\chi(M)$ for II$_1$ factors
satisfying these assumptions.

\heading 2. Some generalities on spectral gap
\endheading

\noindent {\it 2.1. Definition}  ([P06a]). Let $M$ be a II$_1$
factor and $Q\subset M$ a von Neumann subalgebra. We say that $Q$
has {\it spectral gap} in $M$ if $\forall \varepsilon > 0$, $\exists
u_1, ..., u_n\in \Cal U(Q)$ and $\delta > 0$ such that if $x\in M$
satisfies $\|[x,u_i]\|_2 \leq \delta \|x\|_2$ , $\forall i$, then
$\|x-E_{Q'\cap M}(x)\|_2 \leq \varepsilon \|x\|_2$. Note that this
is equivalent to the condition that the representation of $\Cal
U(Q)$ on the Hilbert space $L^2(M\ominus Q'\cap M)=L^2(M)\ominus
L^2(Q'\cap M)$, given by $u \mapsto \text{\rm Ad}u$, has spectral
gap, i.e. it does not weakly contain the trivial representation of
$\Cal U(Q)$.

\vskip .05in

\noindent {\it 2.2. Remarks}. $1^\circ$ Note that if $Q$ is a II$_1$
factor, 2.1 automatically implies $Q$ is a non-Gamma II$_1$ factor.

\vskip .05in

$2^\circ$  A weaker version of spectral gap, is the following:
$\forall \varepsilon > 0$, $\exists$ $u_1, ..., u_n \in \Cal U(Q)$
and $\delta> 0$ such that if $x\in (M)_1$ satisfies $\|[x,u_i]\|_2
\leq \delta $, $\forall i$, then $\|x-E_{Q'\cap M}(x)\|_2 \leq
\varepsilon$.  This latter condition, to which we refer as ``$Q$
{\it has w-spectral gap in} $M$'', is trivially equivalent to the
condition $Q'\cap M^\omega=(Q'\cap M)^\omega$. In turn, spectral gap
is equivalent to $Q'\cap M^2(M)^\omega = L^2(Q'\cap M)^\omega$,
where for a Hilbert space $\Cal H$ and a free ultrafilter $\omega$
on $\Bbb N$, $\Cal H^\omega$ denotes the $\omega$-ultrapower Hilbert
space. It is the stronger form of spectral gap that one usually
checks in concrete examples, but the weaker version is what one
typically needs in deformation-rigidity arguments. As also pointed out 
in (3.4 in [P06]), the difference
between w-spectral gap and spectral gap is in the same vein as the
difference between strong ergodicity and spectral gap of a measure
preserving action of a group on a probability space $\Gamma
\curvearrowright X$, emphasized by K. Schmidt in ([S81]). It is also 
analogue to the
difference between the non inner amenability of a group $\Gamma$ (as
defined by Effros in [E78]) and the property non-Gamma of its
associated II$_1$ factor $L(\Gamma)$, discovered recently by Vaes
([Va09]).

However, in many interesting cases, spectral gap is equivalent to
w-spectral gap, notably in the case $Q\subset Q\overline{\otimes}
S=M$. Indeed, by a result of Connes in [C76], both conditions are
equivalent to $Q$ being non-Gamma. More generally, if $Q_0\subset Q$
is an inclusion of non-Gamma factors with finite index, then
$Q_0\subset Q\subset Q\overline{\otimes} S$ has spectral gap.

But one can easily ``cook up'' examples where they are not. Thus,
let $\theta_k$ be the automorphism of $L(\Bbb F_\infty)$ that takes
the $k$'th generator to its negative and note that $\theta_k
\rightarrow id$. For each $n\geq 1$, let $N_n\subset M_2(\Bbb
C)\otimes L(\Bbb F_\infty)$ be the ``diagonal'' subfactor given by
$\{ x e_{11} + \theta_n(x) e_{22} \mid x\in L(\Bbb F_\infty) \}$,
where $\{e_{ij}\}_{1 \leq i,j \leq 2}$ are the matrix units in
$M_2(\Bbb C) = M_2(\Bbb C)\otimes 1$. Then let $M\simeq L(\Bbb
F_\infty)$ and choose a partition of 1 with projections in $M$,
$\{p_n\}_{n \geq 0}\in M$. For each $n$ we identify $p_nMp_n$ with
$M_2(L(\Bbb F_\infty))$, via some isomorphism $\sigma_n$, by using
Voiculescu's theorem ([V88]), and let $N\subset M$ be the
``diagonal'' subfactor $\{\Sigma_n \sigma_n^{-1}( x e_{11} +
\theta_n(x)e_{22}) \mid x\in L(\Bbb F_\infty) \}$. Then $\xi_m =
\sigma_m^{-1}(e_{12})/\tau(p_m)^{1/2}\in L^2(M)$ satisfy
$\|\xi_m\|_2^2=1/2$ and $\|[\xi_m, x]\|_2 \rightarrow 0$, $\forall
x\in N$. On the other hand, it is trivial to see that if $N$ is a
non-Gamma II$_1$ factor and $N\subset M$ is so that $N'\cap M$ is
atomic with $pMp=Np$ for any minimal projection $p\in N'\cap M$,
then $N$ has w-spectral gap in $M$.

\vskip .05in

$3^\circ$ Several examples of subfactors $Q\subset M$ with spectral
gap are emphasized in [P09]. For instance, if $N_{-1}\subset N_0$ is
an inclusion of II$_1$ factors with finite index, $N_n \nearrow
N_\infty$ is its associated Jones tower and enveloping algebra, and
we assume that either $N_0$ has the property (T), or $N_0$ is merely
non-Gamma but $N_{-1}\subset N_0$ has finite depth, then $N_n$ has
spectral gap in $N_\infty$, $\forall n$ (see $2^\circ$ and $3^\circ$
in Proposition 3.2 on [P09]). Also, if $\Cal G$ is a standard
$\lambda$-lattice ([P94]), $Q$ is a $\text{\rm II}_1$ factor and one
denotes $N_{-1}=M_{-1}^{\Cal G}(Q)\subset M_0^{\Cal G}(Q)=N_0$ the
inclusion of $\text{\rm II}_1$ factors with standard invariant $\Cal
G$, as constructed in ${\text{\rm [P94], [P98]}}$, with $N_n
\nearrow N_\infty$ the associated tower and enveloping algebra, then
$N_n$ has spectral gap in $N_\infty$, $\forall n$ (see 3.2.4$^\circ$
in [P09]).

Another class of examples is obtained below:

\proclaim{2.3. Proposition} Let $N$ be a non-Gamma $\text{\rm II}_1$
factor and $T=N\overline{\otimes} R$. Let $\Gamma
\curvearrowright^\sigma T$ be a cocycle action of a countable group
$\Gamma$ and denote $M=N\overline{\otimes} R \rtimes \Gamma$.

\vskip .05in

$1^\circ$ For any $g$, there exists $t>0$ such that, modulo
perturbation by an inner automorphism of $T$, one has
$\sigma(g)(N)=N^t$, $\sigma(g)(R)=R^{1/t}$, where
$N^t\overline{\otimes} R^{1/t}$ is the unique decomposition of
$N\overline{\otimes} R$ in $\text{\rm [OP03]}$.

\vskip .05in $2^\circ$ Denote $\theta(g)=\sigma(g))_{|N}$ $($cf part
$1^\circ)$. Then $\theta(g)$ is outer, $\forall g\neq e$, iff
$\theta(\Gamma) \cap \overline{\text{\rm Int}(T)}=1$ and iff $N'\cap
M=R$.

\vskip .05in $3^\circ$ Assume the equivalent conditions in $2^\circ$
are satisfied. Then $\rho(g)=\sigma(g)_{|R}$ is outer, $\forall
g\neq e$, iff $N$ satisfies the bicommutant condition $(N'\cap
M)'\cap M=N$ and iff $\sigma$ is centrally free.

\vskip .05in $4^\circ$ Assume the equivalent conditions in $2^\circ$
are satisfied. Then $N$ has spectral gap in $M$ iff $\sigma(\Gamma)$
is closed in $\text{\rm Out}(T)$ and iff $\theta(\Gamma)$ is closed
in $\text{\rm Out}(N^\infty)$.
\endproclaim
\noindent {\it Proof}. Part is a consequence of (Theorem 5.1 in
[P06b]). We leave the proof of $2^\circ-4^\circ$  as an exercise.
\hfill $\square$

\proclaim{2.4. Lemma} Let $Q\subset M$ be an inclusion of factors
and assume that $Q$ has spectral gap in $M$. If $N_n\subset M$ are
von Neumann subalgebras such that $\lim_n \|E_{N_n}(x)-x\|_2=0$,
$\forall x\in M$, then for any $\varepsilon > 0$, there exists $n$
such that $N_n'\cap M \subset_{\varepsilon} Q'\cap M$. Moreover, for
this $n$ we have $(Q'\cap M)'\cap M \subset_{2\varepsilon} (N_n'\cap
M)'\cap M$.
\endproclaim
\vskip .05in \noindent {\it Proof.} Since $Q$ has spectral gap in
$M$, by Lemma 2.2 there exist $\delta>0$ and $u_1, ..., u_m\in \Cal
U(Q)$ such that if $x\in (M)_1$ satisfies $\|[u_i, x]\|_2 < \delta$,
$\forall i$, then $x\in_\varepsilon Q'\cap M$.

By the hypothesis, there exists $n$ such that
$\|E_{N_n}(u_i)-u_i\|_2 \leq \delta/2$, $\forall i$. Thus, if $x\in
N_n'\cap M$ then

$$
\|[x, u_i]\|_2 \leq 2\|x\| \|E_{N_n}(u_i)-u_i\|_2 \leq \delta,
\forall i,
$$
implying that $x\in_\varepsilon Q'\cap M$. The other approximate
inclusion is now automatic, by the triangle inequality.

\hfill $\square$

\heading 3. Virtual strong rigidity results \endheading

\noindent {\it 3.1. Definitions}. $1^\circ$ Recall from [P97], [P01] that a
von Neumann subalgebra $B$ of a II$_1$ factor $M$ is {\it
quasi-regular} in $M$ if the von Neumann algebra generated by its
{\it quasi-normalizer}, defined as the set $q\Cal N(B)$ of all $x\in
M$ with the property that sp$BxB$ is finitely generated both as a
left and as a right $B$-module, is equal to $M$. Note that this is
equivalent to $B'\cap \langle M, e_B \rangle$ being generated by
projections $p$ with the property that both $p$ and $J_MpJ_M$ are
finite in $\langle M, e_B \rangle$. Moreover, by (1.4.2 in [P01]),
$B$ is in fact quasi-regular in $M$ iff it is {\it discrete} in $M$,
in the sense of [ILP96], i.e. $B'\cap \langle M, e_B \rangle$ is
generated by projections that are finite in $\langle M, e_B
\rangle$.

\vskip .05in
\noindent
$2^\circ$ If $N, Q$ are von Neumann subalgebras of the II$_1$ factor $M$, then
$Q$ is {\it discrete over} $N$ ({\it inside $M$}) if $L^2(M)$ is
generated by $Q-N$ Hilbert-bimodules which are finitely generated
over $N$. Arguing as in the proof of (1.4.2 in [P01]), this is
easily seen to be equivalent to $Q'\cap \langle M, e_N \rangle$
being generated by projections that are finite in $\langle M, e_N
\rangle$.

\vskip .05in
\noindent
$3^\circ$ If $N\subset M$, $Q\subset P$ are quasi-regular von Neumann subalgebras of the
II$_1$ factors $M, P$, then an isomorphism $\theta:M\simeq P$ is a
{\it virtual isomorphism} of the inclusions $N\subset M$ and $Q\subset P$, if
$L^2(P)$ is a direct sum of $\theta(N)-Q$ bimodules of finite index,
in other words if $Q$ is discrete over $\theta(N)$ and $\theta(N)$
is discrete over $Q$, inside $P$. Note that this notion becomes relevant only
in the case $N, Q$ have infinite index in their respective ambient algebras.
Indeed, if $N\subset M$, $Q\subset P$ are subfactors of finite Jones index,
then any isomorphism of $M, P$ satisfies the condition.  In turn, for finite index
inclusions of factors we have the following notion of ``virtual isomorphism'':

\vskip .05in
\noindent
$4^\circ$ Let $N_{-1}\subset N_0$, $Q_{-1}\subset Q_0$ be inclusions of
II$_1$ factors with finite Jones index and $N_\infty, Q_\infty$ their
enveloping algebras. An isomorphism $\Delta: N_\infty \simeq Q_\infty$
with the property that $L^2(Q_\infty)$ is generated by $\theta(N_{-1})-Q_{-1}$
Hilbert bimodules which are finite over both $\theta(N_{-1})$ and $Q_{-1}$,
is called a {\it weak equivalence} of $N_{-1}\subset N_0$, $Q_{-1}\subset Q_0$.

\vskip .05in
\noindent
$5^\circ$ Finally, if $N, Q\subset M$ are von Neumann subalgebras of the
II$_1$ factor $M$, then a (non-zero) Hilbert $N-Q$ bimodule $\Cal
H\subset L^2(M)$, which is finitely generated both over $N$ and over
$Q$, is called a {\it virtual conjugacy of} (corners of) {\it $N, Q$ inside $M$}. If
such a bimodule exists, we say that $N,Q$ are {\it virtually conjugate
inside $M$} and write $N \sim_M Q$. Like with the virtual isomorphism
of inclusions defined in $3^\circ$, this relative version ``doesn't
recognise'' between finite index inclusions: if $Q, N$ are any subfactors
of finite index of a II$_1$ factor $M$, then $N \sim_M Q$.

\proclaim{3.2. Lemma} Let $N, Q$ be quasi-regular subfactors of the
$\text{\rm II}_1$ factor $M$. Then the following conditions are
equivalent:

$(a)$ $N$ is discrete over $Q$ and $Q$ is discrete over $N$.

$(b)$ $L^2(M)$ is generated by $N-Q$ bimodules that are finitely
generated both over $N$ and over $Q$.

$(c)$ $N\sim_M Q$.
\endproclaim
\noindent {\it Proof}. Let $\tilde{M}=M_2(M)$ and $\tilde{N}=N\oplus
Q \subset \tilde{M}$. It is immediate to see that condition $(b)$ is
equivalent to $\tilde{N}$ being quasi-regular in $\tilde{M}$, while
the condition $(a)$ is equivalent to $\tilde{N}$ being discrete in
$\tilde{M}$. Thus, $(a) \Leftrightarrow (b)$ by (part $(iii)$ of
1.4.2 in [P01]). On the other hand, we trivially have $(b)
\Rightarrow (c)$ and for the converse, all we need to prove is that
$\tilde{N}$ is quasi-regular in $\tilde{M}$. But since the
quasi-normalizers of $Qe_{11}$ and $Ne_{22}$ in $\tilde{M}$ are
included into the quasi normalizer of $\tilde{N}$ and they generate
the factors $e_{11}Me_{11}, e_{22}Me_{22}$ respectively, and since
$N\sim_M Q$ means there exists a non-zero element $x =
e_{11}xe_{22}\in q\Cal N_{\tilde{M}}(\tilde{N})$, it follows that
the von Neumann algebra generated by $q\Cal
N_{\tilde{M}}(\tilde{N})$ must be all $\tilde{M}$.

\hfill $\square$

\proclaim{3.3. Lemma} Let $N=N_0\subset N_1 \subset N_2 \subset ...
\nearrow M$ be an inductive limit of $\text{\rm II}_1$ factors such
that $[N_{n+1}: N_n]< \infty$. Then $N$ is quasi-regular in $M$. If
in addition $N_n'\cap M$ is a factor, $\forall n\geq 0$, then
$N_{n+1}'\cap M \subset N_{n}'\cap M$ has finite index, $N_n \vee
N_n'\cap M$ is quasi-regular in $M$ and $N_n$ is quasi-regular in
$(N_n'\cap M)'\cap M$, $\forall n\geq 0$.
\endproclaim
\noindent {\it Proof}. Since $e_{N_n}\in \langle M, e_N \rangle$ are
finite projections and tend to 1, the first part follows trivially
from (1.4.2 in [P01]).

To prove the second part, it is clearly sufficient to consider the
case $n=0$. The fact that $N_1'\cap M$ has finite index in $N'\cap
M$ was shown in [P92], but we include a full proof for the sake of
completeness. Thus, denote by $e\in N_1$ the Jones projection for
$N\subset N_1$ and let $P=\{e\}'\cap N\subset N$ be the
corresponding downward basic construction (cf [J82]; see also
[PP83]). Let also $\{m_j\}_j$ be an orthonormal basis of $N$ over
$P$ (cf. [PP83]). Thus, $\sigma_j m_j e m_j^*=1$ and $p_j =
E_N(m_j^*m_j)\in \Cal P(P)$, $\forall j$. Note that $N'\cap M \ni
x\mapsto \Cal E(x)=\Sigma_j m_je x em_j^*$ defines a conditional
expectation of $N'\cap M$ onto $N_1'\cap M$, satisfying $exe=\Cal
E(x)e$ and
$$
\tau(\Cal E(x)y)= \Sigma_j \tau(x em_j^*m_j e y)=\Sigma_j \tau(p_jx
ey) \tag 3.3.1 $$ $$=\Sigma_j \tau(E_{P'\cap M}(p_j)x ey) =[N_1:N]
\tau(xey)=\tau(bxy), \forall y\in N_1'\cap M,
$$
where $b=[N_1:N] E_{N'\cap N_1}(e)$. Notice that $[N_1:N]^{-1}1 \leq
b \leq [N_1:N]1$ and denote $E(x)=\Cal E(b^{-1/2}xb^{-1/2})$. By
$(3.3.1)$, it follows that $E$ defines the trace preserving
conditional expectation of $N'\cap M$ onto $N_1'\cap M$. Moreover,
$e_0=b^{-1/2}eb^{-1/2}$ is a projection and $e_0xe_0=E(x)e_0,
\forall x\in N'\cap M$. Thus, $[N'\cap M:N_1'\cap M] \leq
\tau(e_0)^{-1}=\tau(e)^{-1}=[N_1:N]$.

The fact that $N\vee N'\cap M$ is quasi-regular in $M$ follows as in
[P98], by noticing that if $e_n \in N_n$ is a Jones projection for
$N\subset N_n$, then $\Cal H_n=L^2(N\vee N'\cap M e_n N\vee N'\cap
M)$ is finitely generated both as left and as right Hilbert $N\vee
N'\cap M$-module. Since $\Cal H_n L^2(M)$, we are done.

Finally, to show that $N$ is quasi-regular in $\tilde{N}=(N'\cap
M)'\cap M$, note that one has the commuting square

$$
\CD
N\vee N'\cap M\ @.\subset\ @.M\\
\noalign{\vskip-6pt}
\cup\ @.\ @.\cup\\
\noalign{\vskip-6pt} N\ @.\subset\ @.\tilde{N}
\endCD
$$
Thus, if we denote $e=e^M_{N\vee N'\cap M}$, then $e$ implements the
trace preserving conditional expectation of $\tilde{N}$ onto $N$ and
if $f = \vee \{ueu^* \mid u\in \Cal U(\tilde{N}) \}$, then
$Nf\subset \tilde{N}f \subset^{e} \Cal N=\text{\rm vN}(\tilde{N},
e)$ is the basic construction for $N\subset \tilde{N}$. Since for
any projection $p\in (N\vee N'\cap M)'\cap \langle M, e \rangle$ of
finite trace, $fpf \in f \langle M, e \rangle f$ has finite trace,
commutes with $N$ and its $Tr$-preserving expectation onto $\Cal N$
still has finite trace and commutes with $N$, it follows that
$N'\cap \Cal N$ is generated by finite projections. Thus, $N\subset
\tilde{N}$ is quasi-regular. \hfill $\square$

\vskip .05in

\noindent {\it 3.4. Definition}. An inclusion of II$_1$ factors
$N\subset M$ is {\it proper}, if there exists an increasing sequence
of factors $N \subset N_1 \subset N_2 \subset ... \subset M$ such
that:

\vskip .05in \noindent $(3.4.1)$ $\overline{\cup_n N_n}=M$;

\vskip .05in \noindent $(3.4.2)$ $[N_n:N]< \infty$, $\forall n$;

\vskip .05in \noindent $(3.4.3)$ $N_n$ has spectral gap in $M$,
$\forall n$;

\vskip .05in \noindent $(3.4.4)$ $(N_n'\cap M)'\cap M=N_n$, $\forall
n$.

\proclaim{3.5. Theorem} Let $N\subset M$, $Q\subset P$ be proper
inclusions of factors. Any $\theta: M \simeq P$ implements a virtual isomorphism
between the inclusions $N\subset M$, $Q\subset P$ and
between the inclusions $N\vee N'\cap M \subset M$, $Q\vee Q'\cap P
\subset P$.
\endproclaim
\noindent {\it Proof}. We identify $P$ with $M$, via $\theta$. Let
$N\subset N_1 \subset ... \nearrow M$ and $Q \subset Q_1 \subset ...
\nearrow P$ be inductive limits of subfactors satisfying conditions
3.4. By Lemma 2.4, the spectral gap and bicommutant properties for
$N_n$ and $Q_m$ in $M$, given any $\varepsilon > 0$ there exist $m,
n$ such that $N_m'\cap M \subset_\varepsilon Q_n'\cap
M\subset_{\varepsilon} N'\cap M$ and $N\subset_{2\varepsilon} Q_n
\subset_{2\varepsilon} N_m$. Taking $\varepsilon < 1/4$, it follows
that if we average the projection $e_{Q_n \vee (Q_n'\cap M)}\in
\langle M, e_{Q_n \vee (Q_n'\cap M)} \rangle$ by unitaries in $\Cal
U(N) \Cal U(N_m'\cap M)$ gives a finite non-zero projection. Thus,
by [P03] we get $N\vee (N_m'\cap M) \prec_M Q_n \vee (Q_n'\cap M)$.
Since $N \vee (N_m'\cap M)$ has finite index in $N\vee (N'\cap M)$,
this implies $N\vee (N' \cap M) \prec_M Q_n\vee (Q_n'\cap M)$. But
$Q\vee (Q'\cap M)$ and $Q_n\vee (Q_n'\cap M)$ have finite index over
their common subfactor $Q \vee (Q_n'\cap M)$, thus $N\vee (N' \cap
M) \prec_M Q_n\vee (Q_n'\cap M)$ implies $N\vee (N' \cap M) \prec_M
Q\vee (Q_n'\cap M)$ which in turn implies $N\vee (N' \cap M) \prec_M
Q\vee (Q'\cap M)$, with the latter being regular in $M$. Similarly,
$Q\vee (Q'\cap M) \prec_M N\vee (N' \cap M)$ and the latter is
regular in $M$. By the proof of (Lemma 2 in [P04]), it then follows
that $L^2(M)$ is a direct sum of finite dimensional  $N\vee (N' \cap
M)-Q\vee (Q'\cap M)$ bimodules.

The proof that $_N L^2(M)_Q$ is generated by finite $N-Q$ bimodules
is similar and is left as an exercise.
\hfill $\square$

\proclaim{3.6. Corollary} Let $N_{-1}\subset N_0$, $Q_{-1}\subset
Q_0$ be non-Gamma subfactors with finite depth. Then any isomorphism
of their enveloping algebras $\Delta: N_\infty \simeq Q_\infty$
implements a weak equivalence of $N_{-1}\subset N_0$, $Q_{-1}\subset
Q_0$.
\endproclaim
\noindent {\it Proof.} By [P09], both $N_0\subset N_\infty$ and $Q_0\subset Q_\infty$ are
proper inclusions, so Theorem 3.5 applies.
\hfill $\square$

Let us mention that a ``virtual strong rigidity result'' similar to
Theorem 3.5 was
already shown in [P04]:

\proclaim{3.7. Theorem} Let $N\subset P$, $Q\subset P$ be
irreducible, quasi-regular, rigid inclusions of $\text{\rm II}_1$
factors such that $M$ $($rep. $P)$ has Haagerup property relative to
$N$ $($resp $Q)$. Then any $\theta: M \simeq P$ implements a virtual
isomorphism between the inclusions $N\subset M$, $Q\subset P$.
\endproclaim
\noindent {\it Proof.} The proof of (Lemma 2 and the
Theorem in [P04]) actually shows this. \hfill $\square$

\heading 4. Proper inclusions coming from proper actions \endheading

We have already pointed out in Proposition 2.3 that if a group
$\Gamma$ acts on a s-McDuff factor $T=N\overline{\otimes} R$ such
that its restrictions to the non-Gamma and hyperfinite parts $N, R$
of $T$ are both outer and its image in Out$(T)$ is closed, then $N$
has spectral gap in $M=T \rtimes \Gamma$ and satisfies the
bicommutant condition. We will consider here the following more
specific case:

\vskip .05in

\noindent {\it 4.1 Definition} Let $\Gamma$ be a discrete group and
$N\overline{\otimes} R$ a s-McDuff factor, with $N$ non-Gamma. A
{\it proper} action of $\Gamma$ on $N\overline{\otimes} R$ is a
free, cocycle action with the property that $\sigma$ normalizes $N$
and $R$, its restrictions to $N$ and $R$ are both properly outer and
$\sigma_{|N}$ is closed in Out$(N)$. Note that, any such $\sigma$ is
of the form $\sigma=\theta\otimes \rho$, where $\theta: \Gamma
\hookrightarrow \text{\rm Out}(P)$, $\rho: \Gamma\hookrightarrow
\text{\rm Out}(R)$ are faithful $\Gamma$-kernels with H$^3(\Gamma,
\Bbb T)$-obstruction satisfying $\text{\rm Ob}(\rho)=
\overline{\text{\rm Ob}(\theta)}$.

\proclaim{4.2. Theorem} If $\Gamma$ is a countable amenable group,
$N\overline{\otimes} R$ an s-McDuff factor and $\Gamma
\curvearrowright^\sigma N\overline{\otimes} R$  is a proper action,
then $N\subset M=N\overline{\otimes} R \rtimes \Gamma$ is a proper
inclusion.
\endproclaim
\noindent {\it Proof.} It is easy to see that it is sufficient to
prove the case $N$ is separable. But then the statement follows
immediately from the following:

\proclaim{4.3. Lemma} Let $N$ be a separable $\text{\rm II}_1$
factor, $\Gamma$ a countable amenable group and $\theta:\Gamma
\hookrightarrow \text{\rm Out}(N)$ a faithful $\Gamma$-kernel.

$1^\circ$ There exists a hyperfinite subfactor $R\subset N$, with
$R'\cap N=\Bbb C$, and a lifting $\theta':\Gamma \rightarrow
\text{\rm Aut}(N)$ such that $\theta'(g)(R)=R$, $\forall g\in
\Gamma$, $\theta'(g)\theta'(h)=\text{\rm Ad}(u_{g,h}) \theta'(gh)$,
for some $u_{g,h}\in \Cal U(R)$, $\forall g,h\in \Gamma$, and
$\theta'_{|R}$ implements a faithful $\Gamma$-kernel in $\text{\rm
Out}(R)$.

$2^\circ$ Let $\rho:\Gamma \hookrightarrow \text{\rm Out}(R)$ be a
faithful $\Gamma$-kernel on the hyperfinite $\text{\rm II}_1$ factor
with $\text{\rm Ob}(\rho)=\overline{\text{\rm Ob}(\theta)}$ and
denote $\sigma=\theta\otimes \rho$ the corresponding cocycle action
of $\Gamma$ on $N\otimes R$. Then there exists an increasing
sequence of subfactors with finite index, $N\subset N_1 \subset
...\nearrow N_\infty=N\otimes R \rtimes_{\theta\otimes \rho}
\Gamma$, such that $(N_n'\cap N_\infty)'\cap N_\infty=N_n$, $\forall
n$.
\endproclaim
\vskip .05in \noindent {\it Proof}. Note that $\theta\otimes
\theta^{\o\p}$ implements a cocycle action on $N\otimes N^{\o\p}$.
Denote $\tilde{M}=N\otimes N^{\o\p} \rtimes_{\theta \otimes
\theta^{\o\p}} \Gamma$ the corresponding crossed product II$_1$
factor. All work will be done in this ``ambient'' algebra. We will
prove simultaneously $1^\circ$ and $2^\circ$.

We first treat the case $\Gamma$ is generated by finitely many
elements $g_0=e, g_1, ..., g_n\in \Gamma$. Note that by perturbing
if necessary each $\theta_g$ by an inner automorphism implemented by
a unitary in $N$, we may assume $\sigma$ fixes a subalgebra
$M_{(n+1)\times (n+1)}(\Bbb C)$ generated by matrix units
$\{e_{ij}\}_{i,j}$. Thus, we can ``split'' $N$ as $N=M_{(n+1)\times
(n+1)}(\Bbb C)\otimes P$, while at the same time realizing $P$ as a
locally trivial ``diagonal'' subfactor of $N$, $P=\{\Sigma_{j=0}^n
\theta(g_j)(x)e_{jj} \mid x\in P\}$.

By [P98], we can identify $N\otimes N^{\o\p} \subset \tilde{M}$ with
the {\it symmetric enveloping inclusion} $N \vee N^{\o\p} \subset
\dis N\bt_{e_P}N^{\o\p}$ of $P\subset N,$ as defined in [P98]. By
[P98] again, one can embed the Jones tower of factors $P\subset N
\subset N_1\subset ...$, with its enveloping algebra $N_\infty$, as
subfactors in the symmetric enveloping algebra, $N\otimes 1 =
N\subset N_1 \subset ... \nearrow N_\infty \subset \tilde{M}$.
Moreover, any such embedding of the tower of factors comes from a
tunnel of factors $N\supset P \supset P_1 \supset ... $, by taking
$N_n={P_{n-1}^{\o\p}}'\cap \tilde{M},$ where $^{\o\p}$ is the
canonical symmetry of $\tilde{M}$, which leaves $e_P$ fixed and
sends $N$ onto $N^{\o\p}$ (see [P98]).

Other observation concerning the symmetric enveloping algebra
$\tilde{M}$ that will be useful are the following (cf. [P98]):
$N'\cap \tilde{M}=N^{\o\p}$ and $(N^{\o\p})'\cap \tilde{M}=N$; more
generally $(N_n'\cap \tilde{M})'\cap \tilde{M}=N_n$, $\forall n\in
\Bbb Z$, where $N_0=N$, $N_{-1}=P$ and $N_k=P_{-k-1}$ for $k\leq
-2$; given any other locally trivial ``diagonal'' subfactor
$P^0\subset N$, corresponding to another set of generators (with
possible multiplicities) $h_0=e, h_1, ..., h_m\in \Gamma$ and
appropriate inner perturbations of $\theta(h_i)$, one can
alternatively view $\tilde{M}$ as the symmetric enveloping algebra
$\dis N\bt_{e_{P^0}}N^{\o\p}$; given any ``level'' $k\in \Bbb Z$ of
the tower-tunnel $N_k$, one can identify $N\vee N^{\o\p} \subset
\dis N\bt_{e_P}N^{\o\p}= N\otimes N^{\o\p} \rtimes_{\theta\otimes
\theta^{\o\p}} \Gamma$ with $N^t \otimes {N^{1/t}}^{\o\p}\subset N^t
\otimes {N^{1/t}}^{\o\p} \rtimes_{\theta^t\otimes
{\theta^{1/t}}^{\o\p}} \Gamma$, which we view also as $N_k \vee
N_k'\cap \tilde{M}\subset \tilde{M}$ and as the amplification by
$1/t^2$ of the symmetric enveloping inclusion $N_k \vee N_k^{\o\p}
\subset \dis N_k\bt_{e_{N_{k-1}}}N_k^{\o\p}$, where $t=(n+1)^k$.

For clarity, we'll first give a short argument in the case $\Gamma$
is strongly amenable with respect to the set of generators $g_0,
g_1, ..., g_n$ (in the sense of [P89a]). Thus, by (proof of 2.1 in
[P89b]), the tunnel $N\supset P \supset P_{1} \supset P_{2} ...$ can
be chosen in this case such that if we denote $R=\overline{\cup_n
P_{n}'\cap N}$, $S_k=\overline{\cup_n P_{n}'\cap P_k}$ then, $R'\cap
\tilde{M} = N'\cap \tilde{M}=N^{\o\p}$ and $S_k'\cap \tilde{M} =
P_k'\cap \tilde{M}$, $k\geq 0$. Equivalently, this means that after
perturbing if necessary each $\theta(g_j)$ by an inner automorphism
of $N$, there exists a $\theta$-invariant hyperfinite subfactor $Q$
of $N$, with relative commutant in $\tilde{M}$ equal to $N^{\o\p}$,
such that $\theta(g)_{|Q}\in \text{\rm Aut}(Q)$, $g\in \Gamma$,
implements a $\Gamma$-kernel, $\Gamma \hookrightarrow \text{\rm
Out}(Q)$, and such that $R = M_{(n+1) \times (n+1)}(\Bbb C)\otimes
Q$, $S=\{\Sigma_j \theta(g_j)(x)e_{jj} \mid x\in Q\}$ (see proof of
2.1 in [P89b]). Keeping in mind that if we let
$N_n={P_{n-1}^{\o\p}}'\cap \tilde{M},$ $N_\infty = \overline{\cup_n
N_n}\subset \tilde{M}$ the representations inside $\tilde{M}$ of the
tower and enveloping algebra of $P\subset N$ corresponding to the
choice of tunnel $P_n$, then $(N_1'\cap N_\infty \subset N'\cap
N_\infty)$ coincides with $(S\subset R)^{\o\p}$ (see [P92] or
[P98]). Also, if we denote by $\rho: \Gamma \hookrightarrow
\text{\rm Out}(R^{\o\p})$ the corresponding $\Gamma$-kernel, then by
(Sec 3 in [P98]), it follows that $\theta \otimes \theta^{\o\p}$
implements a cocycle action of $\Gamma$ on $N \vee N^{\o\p}$ which
leaves $N \vee R^{\o\p}$ invariant, and such that $N \vee (N'\cap
N_\infty) \subset N_\infty$ is isomorphic $N \otimes R^{\o\p}
\rtimes_{\theta \otimes \rho} \Gamma$. Moreover, from the way we
chose the tunnel, we have
$$
(N_n'\cap N_\infty)'\cap N_\infty=(S_{n-1}^{\o\p})'\cap  N_\infty
\subset (S_{n-1}^{\o\p})'\cap \tilde{M} $$ $$ =(S_{n-1}'\cap
\tilde{M})^{\o\p} = (P_{n-1}'\cap \tilde{M})^{\o\p} = N_n,
$$
which together with the trivial inclusion $(N_n'\cap N_\infty)'\cap
N_\infty \supset N_n$ gives $(N_n'\cap N_\infty)'\cap N_\infty=N_n$.

The proof of the lemma in its full generality is very similar, but
instead of a Jones tunnel of factors $P_n$, we now use (the proof of
7.1 and Remark 7.2.1$^\circ$ in [P98]) in combination with (the
proof of 2.1 in [P89b]), to construct recursively a decreasing
sequence of subfactors $P_n$, with each $P_n$ obtained as a downward
basic construction of a suitable local inclusion $P_{n-1}p \subset
pN_m^0p$, for some $m$ and $p\in P_{n-1}'\cap N_m^0$, where
$P\subset N\subset N_1^0\subset N_2^0 \subset ... \subset \tilde{M}$
is a representation of the Jones tower inside $\tilde{M}$, which we
choose once for all, from the beginning.

Let us consider first the case $\Gamma$ generated by a finite set
$g_0=e, g_1, ..., g_n$ and let $P\subset N \subset \tilde{M}$ be
defined as before. Let $\{x_n\}_n \subset \tilde{M}$ be dense in the
norm $\|\cdot \|_2$ in the unit ball of $\tilde{M}$. Assume we have
constructed finite index subfactors $N=P_{-1} \supset P=P_0 \supset
P_1 ... \supset P_{n-1}$ such that $E_{(P_j'\cap P_k)'\cap
\tilde{M}}(x_i)\in_{2^{-j}} P_k'\cap \tilde{M}$, $\forall 1\leq i
\leq j \leq n-1$, $-1\leq k \leq j \leq n-1$, and such that each
$P_{j}\subset N$ is a locally trivial subfactor, corresponding to
some finite subset $e\in K_{j}\subset \Gamma$ and multiplicities
$k_{g, j}\geq 1, g\in K_j$.

By (4.10 in [P98]), the inclusion $N\supset P_{n-1}=P^{-1}_0$ admits
a tunnel $N\supset P^{-1}_0 \supset P^{-1}_1 \supset ... \supset
P^{-1}_k$ such that $E_{({P_k^{-1}}'\cap P_{n-1})'\cap
\tilde{M}}(x_j)\in_{2^{-n-1}} {P_k^{-1}}'\cap \tilde{M}$. On the
other hand, since $P_k^{-1} \subset N$ is locally trivial, given by
automorphisms in $\Gamma$, there exists a factor $N\subset N_1^{-1}
\subset \tilde{M}$ such that  $(P_k^{-1})^{\o\p} \subset N^{\o\p}
\subset {P_k^{-1}}'\cap \tilde{M}$ is a basic construction (from
properties mentioned before), let $e\in {P_k^{-1}}'\cap
{(P^{-1}_k)^{\o\p}}'$ be the corresponding Jones projection. Thus,
${P^{-1}_k}'\cap \tilde{M} = \text{\rm sp} N^{\o\p} e N^{\o\p}$. By
(6.1 in [P98]), for any $\delta>0$ there exist $m$, a projection
$p\in {P^{-1}_k}'\cap N_m$ and a downward basic construction $P_n
\subset P^{-1}_k \simeq P^{-1}_kp \subset pN_mp$ such that
$\|E_{({P_n}'\cap P_k)'\cap \tilde{M}}(e)-E_{P_k'\cap
\tilde{M}}(e)\|_2 \leq \delta$. Since

$$
\|E_{({P^{-1}_k}'\cap P_j)'\cap \tilde{M}}(x_i)-E_{{P^{-1}_k}'\cap
\tilde{M}}(x_i)\|_2 \leq 2^{-n-1}
$$
and ${P^{-1}_k}'\cap \tilde{M} =\text{\rm sp} N^{\o\p} e N^{\o\p}$,
$P_j'\cap \tilde{M}=\text{\rm sp} N^{\o\p} E_{P_j'\cap \tilde{M}}(e)
N^{\o\p}$, it follows that for $n-1 \geq j\geq -1$, we have

$$
\|E_{({P_n}'\cap P_j)'\cap \tilde{M}}(x_i)-E_{P_j'\cap
\tilde{M}}(x_i)\|_2
$$
$$
\leq \|E_{({P_n}'\cap P-j)'\cap \tilde{M}}(x_i) - E_{({P_n}'\cap
N)'\cap ({P^{-1}_k}'\cap \tilde{M})}(x_i)\|_2 $$ $$
+\|E_{({P_n}'\cap P_j)'\cap \tilde{M}}(E_{{P^{-1}_k}'\cap
\tilde{M}}(x_i)-E_{P_j'\cap \tilde{M}}(x_i)\|_2
$$
$$
\leq 2^{-n-1} + \|E_{({P_n}'\cap P_j)'\cap
\tilde{M}}(E_{{P^{-1}_k}'\cap \tilde{M}}(x_i))-E_{P_j'\cap
\tilde{M}}(x_i)\|_2,
$$
which for small enough $\delta$ gives $\|E_{({P_n}'\cap P_j)'\cap
\tilde{M}}(x_i)-E_{P_j'\cap \tilde{M}}(x_i)\|_2 \leq 2^{-n}$,
$\forall i \leq n$.

We now define $N_n={P_n^{\o\p}}' \cap \tilde{M}$,
$N_\infty=\overline{\cup_n N_n}$, $R=\overline{\cup_n P_n'\cap N}$,
$S_j=\overline{\cup_n P_n'\cap P_j}$. Then $N_\infty$ is generated
by $N, e_P, R^{\o\p}$. Moreover, by the way $P_n$ were chosen, we
also have $R'\cap \tilde{M}=N^{\o\p}$, $S_j'\cap \tilde{M}=P_j'\cap
\tilde{M}$ and a non-degenerate commuting square of factors

$$
\CD
P\ @.\subset\ @.N\\
\noalign{\vskip-6pt}
\cup\ @.\ @.\cup\\
\noalign{\vskip-6pt} S\ @.\subset\ @.R
\endCD
$$

As in the proof of (2.1 in [P89b]), this implies that there exist
unitary elements $u_i\in N$ such that $\text{\rm
Ad}(u_i)\theta_{g_i}$ normalize $R$. Thus, after some inner
perturbations, $\theta$ normalizes $R$. Moreover, the condition
$R'\cap \tilde{M}=N^{\o\p}$ implies that the restriction of $\theta$
to $R$ is a faithful $\Gamma$-kernel (i.e. $\theta(g)_{|R}$ is inner
iff $g=e$). Altogether, if we denote by $\rho$ the restriction of
$\theta^{\o\p}$ to $R^{\o\p}$, then $\theta\otimes \theta^{\o\p}$
normalizes $N_\infty$, the inclusion $N\vee (N'\cap N_\infty)
\subset N_\infty$ is isomorphic to $N \otimes R^{\o\p}
\rtimes_{\theta\otimes \rho} \Gamma$ and if $S_n=P_n\cap
S=\overline{\cup_k P_k'\cap P_n}$, then
$$
(N_n'\cap N_\infty)'\cap N_\infty=(S_{n-1}^{\o\p})'\cap N_\infty =
$$
$$
=({P_{n-1}^{\o\p}}'\cap \tilde{M})\cap  N_\infty = (P_{n-1}'\cap
\tilde{M})^{\o\p} \cap N_\infty = N_n,
$$
thus proving the statement in the case $\Gamma$ is finitely
generated.

Finally, let's settle the case $\Gamma=\{g_n\}_n$ is infinitely
generated. Let $\Gamma_n\subset \Gamma$ be the subgroup generated by
$g_0=e, g_1, ..., g_m$ and denote $\tilde{M}_m=N\vee N^{\o\p}
\rtimes_{\sigma} \Gamma_m \subset \tilde{M}$. We then construct
$N\supset P \supset P_1 ...$ recursively, so that to satisfy
$E_{(P_j'\cap P_i)'\cap \tilde{M}_j}(x_i)\in_{2^{-j}} P_i'\cap
\tilde{M}$, $\forall 1\leq i \leq j \leq n$. If we let $N_n,
N_\infty$, $R, S_n$ be defined as before, then the statement follows
in its full generality.

\hfill $\square$

Let us also mention that by combining Lemma 4.3 with results in
[Oc83], one obtains a generalization of the vanishing 2-cohomology
(Theorem 2.1 in [P89b]):

\proclaim{4.4. Corollary} Let $N$ be a $\text{\rm II}_1$ factor,
$\Gamma$ a countable amenable group and $\theta:\Gamma
\hookrightarrow \text{\rm Out}(N)$ a faithful $\Gamma$-kernel with
trivial obstruction, $\text{\rm Ob}(\theta)=1$. Then there exists a
lifting $\theta'(g)\in \text{\rm Aut}(N)$, $g\in \Gamma$, such that
$\theta'(g)\theta'(h)=\theta'(gh)$, $\forall g,h\in \Gamma$.
\endproclaim
\noindent {\it Proof}. It is clearly sufficient to prove the result
in the case $N$ is separable. By 5.2.$1^\circ$, there exists a
hyperfinite subfactor $R\subset N$ and a lifting of the
$\Gamma$-kernel $\theta$ to automorphisms $\theta(g)$ that leave $R$
invarian and  satisfy $\theta(g)\theta(h)=\text{\rm Ad}(u_{g,h})
\theta(gh)$, for some $u_{g,h}\in \Cal U(R)$, $\forall g,h\in
\Gamma$. By [Oc83], there exist $w_g\in \Cal U(R)$ such that
$u_{g,h}=\theta(g)(w_h^*)w_g^*w_{gh}$, modulo scalars, $\forall
g,h\in \Gamma$. But then, $\theta'(g)=\text{\rm Ad}w_g \circ
\theta(g)$ satisfy the desired conditions.

\hfill $\square$

\heading 5. Strong rigidity for factors arising from proper
actions\endheading

We now apply the results in the previous sections to derive a
strong rigidity result for crossed product II$_1$ factors arising
from proper actions.

\proclaim{5.1. Theorem} Let $\Gamma_i$ be countable amenable group,
$\sigma_i: \Gamma_i \rightarrow \text{\rm Aut}(N_i
\overline{\otimes} R_i)$ proper cocycle action on the s-McDuff
factor $N_i \overline{\otimes} R_i$, leaving the non-Gamma factor
$N_i$ invariant, and denote $M_i=N_i\overline{\otimes} R_i
\rtimes_{\sigma_i} \Gamma_i$ the corresponding crossed product
factor, $i=0,1$. Let $\Delta: M_0 \simeq M_1$. Then $\Delta$
implements a virtual cocycle conjugacy of $\sigma_0, \sigma_1$, i.e.
a virtual isomorphism of the inclusions $N_0\overline{\otimes}R_0
\subset M_0$, $N_1\overline{\otimes} R_1 \subset M_1$.

If in addition $\Gamma_0, \Gamma_1$ are torsion free, then $\Delta$
implements a stable cocycle conjugacy of $\sigma_0, \sigma_1$. More
precisely, there exists $u\in \Cal U(M_1)$, $\{v_g\}_g \subset \Cal
U(N_1\overline{\otimes} R_1),$ an isomorphism $\delta: \Gamma_0 \simeq
\Gamma_1$ and a splitting $N_1^t\overline{\otimes} R_1^{1/t}$ of $N_1 \overline{\otimes}
R_1$, for some $t> 0$, such that if we denote $\Delta'=\text{\rm
Ad}u \circ \Delta$ and by $\{u^i_g\}_g$ the canonical unitaries in
$M_i=N_i\overline{\otimes} R_i \rtimes \Gamma_i$, then we have:

\vskip .05in \noindent $(a)$ $\Delta'(N_0)=N_1^t$,
$\Delta'(R_0)=R_1^{1/t}$;

\vskip .05in \noindent $(b)$ $\Delta'(x_g u^0_g)= \Delta(x_g)
v_{\delta(g)}u^1_{\delta(g)}$, $\forall x_g\in N_0\overline{\otimes} R_0$,
where $\Delta: N_0\overline{\otimes} R_0 \simeq N_1\overline{\otimes} R_1
=N_1^t \overline{\otimes}
R_1^{1/t}$ is the isomorphism implemented by the restriction of
$\Delta'$ to $N_0\overline{\otimes} R_0$.

\vskip .05in In particular, $\Delta$ implements an outer conjugacy
between the kernel $(\Gamma_0, \theta_0)$ and $(\Gamma_1,
\theta_1^t)$.
\endproclaim
\vskip .05in \noindent {\it Proof.} By Theorem 4.2, the inclusions
$N_i\subset M_i$ are proper, so by Theorem 3.5, $\Delta$ implements
a virtual isomorphism of $N_0\overline{\otimes} R_0\subset M_0$ and
$N_1\overline{\otimes} R_1\subset M_1$. If in additions $\Gamma_i$
are torsion free, then by (Lemma 2 in [P04]) it follows that
$N_0\overline{\otimes} R_0=N_0\vee (N_0' \cap M)$ and
$N_1\overline{\otimes} R_1=N_1\vee (N_1'\cap M)$ are conjugate by a
unitary in $M$. The rest of the statement follows from (6.5.1 in
[P06]). \hfill $\square$

\vskip .05in \noindent {\bf 5.2. Remark}. Let $\Gamma$ be strongly
amenable with respect to the set of generators $g_0=1, g_1, ..., g_n
\in \Gamma$ and let $\alpha \in \text{\rm H}^3(\Gamma, \Bbb T)$. If
$Q$ is a II$_1$ factor, then the proof of Lemma 4.2 shows that the
universal construction of subfactors in ([P94], [P98]) gives rise to
a canonical faithful $\Gamma$-kernel $\theta: \Gamma \hookrightarrow
\text{\rm Out}(N)$ with Ob$(\theta)=\alpha$, where $N=R'\cap (Q
\overline{\otimes} S *_S R)$. Moreover, [PS00] shows that if $Q=L(\Bbb
F_\infty)$ then $N\simeq L(\Bbb F_\infty)$. Thus, the considerations
in 4.2 and 5.1 give canonical faithful $\Gamma$-kernels with given
obstruction on the free group factor $L(\Bbb F_\infty)$, with
discrete closure in Out$(L(\Bbb F_\infty))$. Moreover, by 5.1, two
such $\Gamma$-kernels $\theta:\Gamma \hookrightarrow \text{\rm
Out}(L(\Bbb F_\infty))$ are cocycle conjugate iff the crossed
product factors $L(\Bbb F_\infty) \overline{\otimes} R \rtimes_{\theta \otimes
\rho} \Gamma$ are isomorphic, where $\rho:\Gamma \hookrightarrow
\text{\rm Out}(R)$ is a model $\Gamma$-kernel on the hyperfinite
II$_1$ factor with obstruction Ob$(\rho)$ equal to
$\overline{\alpha}$. Noticing that given any $\Gamma$-kernel
$\theta:\Gamma \hookrightarrow \text{\rm Out}(N)$ we have a natural
identification $(N \overline{\otimes} R \rtimes_{\theta \otimes \rho}
\Gamma)^{\o\p} = N^{\o\p}\overline{\otimes} R^{\o\p} \rtimes_{\theta^{\o\p}
\otimes \rho^{\o\p}} \Gamma$ and that
Ob$(\theta^{\o\p})=\overline{\text{\rm Ob}(\theta)}$, by 5.1 it
follows that if $\alpha\neq \overline{\alpha}$, then the factor
$\tilde{M}=L(\Bbb F_\infty) \overline{\otimes} R \rtimes_{\theta \otimes \rho}
\Gamma$ is not anti-isomorphic to itself. This situation occurs if
for instance $\Gamma=\Bbb Z^3$ and one takes $\alpha\in \text{\rm
H}^3(\Bbb Z^3, \Bbb T)=\Bbb T$, with $\alpha \neq \pm 1$. Similarly,
if $\alpha^n$ are distinct for $n=1,2,...$ (which in this case
amounts to $\alpha = \exp(2\pi i t)$ with $t$ irrational) then the
tensor products $\tilde{M}^{\otimes n}$ are mutually non-isomorphic.
Note that these examples of II$_1$ factors are very similar in
spirit with the ones in [C75].

\vskip .05in

\noindent {\bf 5.3. Remark}. Note that if $P\subset N$ is an
inclusion of non-Gamma factors with finite index then, while the
condition that $N$ has spectral gap in the enveloping algebra
$N_\infty$ is automatic under the finite depth assumption, this is
no longer true in the infinite depth case, not even under the
``strong amenability'' condition on the standard graph $\Cal
G_{P,N}$. This is already clear from 2.3, 4.2, 5.1, where we see
(implicitly) that if $P\subset N$ is a locally trivial subfactor
given by a $\Gamma$-kernel of a finitely generated group, then $N$
has spectral gap in $N_\infty$ iff $\Gamma$ is closed in Out$(N)$.
If $N=L(\Bbb F_n)$, for some $2\leq n \leq \infty$, then given any
finitely generated group $\Gamma$ that can be embedded faithfully
into a unitary group $\Cal U(m)$ (with no scalar values other than
1) it is easy to provide two examples of faithful $\Gamma$-kernels
$\sigma_0, \sigma_1$ on $L(\Bbb F_n)$, one with compact closure in
Out$(L(\Bbb F_n))$, the other one discrete. For instance, one can
take $\sigma_0 : \Gamma \rightarrow \text{\rm Aut}(R)$ to be a
compact product action, obtained from the restriction to $\Gamma$ of
the diagonal action (via Ad) of $O(n)/\Bbb T$ on
$R=\overline{\otimes}_k (M_{m \times m}(\Bbb C))_k$, extended to
$L(\Bbb F_n)=R * L(\Bbb F_{n-1})$ by $\sigma_0=\sigma_0 * 1$.  While
$\sigma_1$ can be taken a Bernoulli action $\Gamma \curvearrowright
R=\overline{\otimes}_{g\in \Gamma} (M_{m \times m}(\Bbb C))_g$,
extended to $L(\Bbb F_n)=R* L(\Bbb F_{n-1})$ by $\sigma_1=\sigma_1
* 1$.

In particular, one can construct two subfactors $P^i\subset N,
i=0,1,$ of the free group factor $N=L(\Bbb F_n)$, both with index 4
and graph $A^{(1)}_\infty=A_{-\infty, \infty}$, the first one with
the property that $N$ has spectral gap in the enveloping factor
$N^0_{\infty}$ of $P^0\subset N$, the second one so that $N$ doesn't
have spectral gap in the enveloping factor $N^1_{\infty}$ of $P^1
\subset N$.

Similarly, by using a construction in ([P89a]), one has two
subfactors $P^i\subset N^i\simeq L(\Bbb F_{2n-1})$, $i=0,1$, of
index $4$ and graph $D_\infty$, such that $N^0$ has spectral gap in
$N^0_{\infty}$ while $N^1$ doesn't have spectral gap in
$N^1_{\infty}$, as follows: Let $\theta_0$ be a properly outer
action of the infinite dihedral group $\Gamma=\Bbb Z/2\Bbb Z * \Bbb
Z/2\Bbb Z$ on $L(\Bbb F_n)$ such that $\theta_0(\Gamma)$ has compact
closure in Out$(L(\Bbb F_n))$. Let $\theta_1$ to be a properly
outer action of $\Gamma$ with discrete closure in  Out$(L(\Bbb
F_n))$. In each case denote by $\sigma_i, \rho_i$ the period 2
automorphisms of $L(\Bbb F_n)$ generating $\theta_i(\Gamma)$. Then
denote $(P^i\subset N^i)=(L(\Bbb F_n)^{\sigma_i} \subset L(\Bbb
F_n)\rtimes \rho_i)$. By [P89a], both subfactors have graph
$D_\infty$. On the other hand, it is easy to check that $N^0$ doesn
not have spectral gap in $N^0_{\infty}$, while $N^1$ does have
spectral gap in $N^1_{\infty}$. Also, by [R92] or [PS00], it follows
that $N^0, N^1 \simeq L(\Bbb F_{2n-1})$.

In conclusion, we see that for non-Gamma subfactors of finite index
$P\subset N$, the condition that $N$ has spectral gap in the
enveloping algebra $N_\infty$ is an isomorphism invariant for
$P\subset N$. The spectral gap property means that the semigroup of
outer symmetries of $N$ (i.e. correspondences, or endomorphisms on
$N^\infty=N\overline{\otimes} \Cal B(\ell^2\Bbb N)$), generated by
the irreducible direct summands of the Hilbert bimodule
$_NL^2(N_\infty)_N$, is closed in End$(N^\infty)/\text{\rm
Int}(N^\infty)$. At the opposite end, we have the case when this
semigroup is precompact.

\vskip .05in

\noindent {\bf 5.4. Remark}. The notion of proper action
on a s-McDuff factor, $\Gamma \curvearrowright^\sigma N\overline{\otimes} R$,
doesn't in fact need the condition that $\sigma$
normalizes $N, R$. Indeed, by 2.3.1$^\circ$ this
is automatic, if one considers instead the cocycle action implemented by $\sigma$
on the II$_\infty$ factors $N^\infty\overline{\otimes} R^\infty$. This, of course,
entails a possible appearence of trace scaling restriction
$\theta=\sigma_{|N^\infty}$,
$\rho=\sigma_{|R^\infty}$, with mod$\theta(g)=
\text{\rm mod}(\rho(g))^{-1}$, $\forall g\in \Gamma$. Theorem 4.3
holds true for these general proper actions
of countable amenable groups on II$_\infty$ s-McDuff factors.
The proof relies on the appropriate generalisation of
Lemma 4.3 to $\Gamma$-kernels in Out$(N^\infty)$, whose
proof is very similar in spirit to the
proof of 4.3, but more ``painful'' to formalise ...
Consequently, Theorem 5.1 holds true in this generality as well.

\vskip .05in

\noindent {\bf 5.5. Remark}. Finally, let us
point out that a generalisation of Lemma 4.3 in the abstract framework of
subfactors holds true, as follows: Given inclusion of separable factors
of type II$_1$, $N_{-1}\subset N_0$, with finite index and amenable graph
(not necessarily extremal!), there exists a decreasing
sequence of subfactor  $N_0\supset N_{-1} \supset N_{-2}\supset ...$
such that if we denote $R_0=\overline{\cup_n N_{-n}'\cap N_0}$,
$R_{-1}=\overline{\cup_n N_{-n}'\cap N_{-1}}$, then $R_i'\cap N_n=N_i'\cap N_n$,
for any $i=0, -1$ and $n\geq 0$, where $N_{-1}\subset N_0\subset N_1 \subset ....$
is the Jones tower. The proof is similar to the proof of 4.3.2$^\circ$. This
can be used to show that if if $N_0$ has spectral gap in  $N_\infty$ then
$N_0\subset N_\infty$ is a proper inclusion, i.e. satisfies conditions 3.3.

\head  References \endhead

\item{[C75]} A. Connes: {\it Sur la classification des facteurs
de type} II, C. R. Acad. Sci. Paris {\bf 281} (1975), 13-15.

\item{[C76]} A. Connes: {\it Classification of injective factors}, Ann.
Math. {\bf 104} (1976),  73-115.

\item{[E75]} E. Effros: {\it Property} $\Gamma$ {\it and inner amenability},
Proc. Amer. Math. Soc. {\bf 47} (1975), 483-486.

\item{[GHJ90]} F. Goodman, P. de la Harpe, V.F.R. Jones: ``Coxeter graphs and
towers of algebras'', MSRI Publications {\bf 14}, Springer, 1989.

\item{[J79]} V.F.R. Jones: {\it A} II$_1$ {\it factor
anti-isomorphic to itself but without involuntary
antiautomorphisms}, Math. Scand. {\bf 46} (1980), 103-117.

\item{[J80]} V.F.R. Jones: Actions of finite groups on the
hyperfinite II$_1$ factor, Mem. A.M.S. {\bf 28}, 1980.

\item{[J82]} V.F.R. Jones: {\it Index for subfactors},
Invent. Math {\bf 72} (1983), 1-25.

\item{[K93]} Y. Kawahigashi: {\it Centrally trivial automorphisms and
an analogue of Connes's} $\chi(M)$ {\it for subfactors}, Duke Math.
J. {\bf 71} (1993), 93-118.

\item{[M70]} D. McDuff: {\it Central sequences and the hyperfinite factor}, Proc.
London Math. Soc. {\bf 21} (1970), 443–461.

\item{[Oc83]} A. Ocneanu: Actions of discrete amenable groups on factors,
Lecture Notes in Math {\bf 1138}, Springer Verlag,
Berlin-Heidelberg-New York, 1985.

\item{[Oc87]} A. Ocneanu: {\it Quantized groups, string algebras and Galois
theory for von Neumann algebras}. In ``Operator Algebras and
Applications'', London Math. Soc. Lect. Notes Series, Vol. {\bf
136}, 1988, pp. 119-172.

\item{[OP03]} N. Ozawa, S. Popa: {\it Some prime factorization results for type}
II$_1$ {\it factors}, Invent. Math., {\bf 156} (2004), 223-234.

\item{[OP07]} N. Ozawa, S. Popa: {\it On a class of $\text{\rm II}_1$
factors with at most one Cartan subalgebra}, math.OA/0706.3623, to
appear in Annals of Mathematics

\item{[PP83]} M. Pimsner, S. Popa: {\it Entropy and index for subfactors}, Ann.
Sc. Ec. Norm. Sup. {\bf 19} (1986), 57-106.

\item{[P89a]} S. Popa: {\it Sur la classification des sous-facteurs
d'indice fini du facteur hyperfini}, Comptes Rendus de l'Acad\'emie
de Sciences de Paris, {\bf 311} (1990), 95-100.

\item{[P89b]} S. Popa: {\it Sousfacteurs, actions des groupes et cohomologie},
Serie I, Comptes Rend. Acad. Sci. Paris, {\bf 309} (1989), 771-776.

\item{[P90]} S. Popa: {\it Markov traces on universal
Jones algebras and subfactors of finite index}, Invent. Math. {\bf
111} (1993), 375-405.

\item{[P91]} S. Popa: {\it Classification of amenable subfactors of type II},
Acta Mathematica, {\bf 172} (1994), 163-255.

\item{[P92]} S. Popa: {\it Classification of actions of
discrete amenable groups on amenable subfactors of type} II, IHES preprint 1992,
to appear in Intern. J. Math., 2009.

\item{[P94]} S. Popa, {\it An axiomatization of the lattice of
higher relative commutants of a subfactor}, Invent. Math., {\bf 120}
(1995), 427-445.

\item{[P97]} S. Popa: {\it Some properties of the symmetric enveloping algebras
with applications to amenability and property} T, Documenta
Mathematica, {\bf 4} (1999), 665-744.

\item{[98]} S. Popa: {\it Universal construction of subfactors}, J.
reine angew. Math., {\bf 543} (2002), 39-81.

\item{[P01]} S. Popa: {\it On a class of type} II$_1$ {\it factors with Betti
numbers invariants}, Ann. of Math {\bf 163} (2006), 809-899.

\item{[P03]} S. Popa: {\it Strong Rigidity of} II$_1$ {\it Factors
Arising from Malleable Actions of $w$-Rigid Groups} I, Invent. Math.
{\bf 165} (2006), 369-408 (math.OA/0305306).

\item{[P04a]} S. Popa: {\it Strong Rigidity of} II$_1$ {\it Factors
Arising from Malleable Actions of $w$-Rigid Groups} II, Invent. Math.
{\bf 165} (2006), 409-453. (math.OA/0407137).

\item{[P04b]} S. Popa: {\it A unique decomposition result for} HT
{\it factors with torsion free core}, J. Fnal. Analysis {\bf 242}
(2007), 519-525 (math.OA/0401138).

\item{[P06a]} S. Popa: {\it On the superrigidity
of malleable actions with spectral gap}, J. Amer. Math. Soc. {\bf
21} (2008), 981-1000 (math.GR/0608429).

\item{[P06b]} S. Popa: {\it On Ozawa's Property for Free Group Factors},
Int. Math. Res. Notices (2007) Vol. {\bf 2007}, article ID rnm036,
10 pages, doi:10.1093/imrn/rnm036 published on June 22, 2007
(math.OA/0608451)

\item{[P06c]} S. Popa: {\it Deformation and rigidity for group actions
and von Neumann algebras}, in ``Proceedings of the International
Congress of Mathematicians'' (Madrid 2006), Volume I, EMS Publishing
House, Zurich 2006/2007, pp. 445-479.

\item{[P09]} S. Popa: {\it On spectral gap rigidity and Connes invariant
$\chi(M)$}, math.OA/0909.5639

\item{[PS00]} S. Popa, D. Shlyakhtenko: {\it Universal properties of}
$L(F_\infty)$ {\it in subfactor theory}, Acta Mathematica, {\bf 191}
(2003), 225-257.

\item{[R92]} F. Radulescu: {\it Random matrices, amalgamated free products and
  subfactors of the von Neumann algebra of a free group, of noninteger
  index}, Invent. Math. {\bf 115} (1994), 347-389.

\item{[S81]} K. Schmidt: {\it Amenability, Kazhdan's property T , strong ergodicity
and invariant means for ergodic group-actions}, Ergodic Theory
Dynamical Systems {\bf 1} (1981), 223-236.

\item{[Va09]} S. Vaes: {\it An inner amenable group whose von Neumann algebra
does not have property Gamma}, math.OA/0909.1485

\item{[V89]} D. Voiculescu: {\it Circular and semicircular systems and free
product factors}, in ``Operator algebras, unitary representations,
enveloping algebras, and invariant theory'' (Paris, 1989), Prog.
Math. {\bf 92}, Birkhäuser, Boston, 1990, 45–60.

\enddocument